\documentclass{amsart}

\usepackage{xypic}
\input xy
\xyoption{all}
\usepackage{epsfig}
\usepackage{amsthm}
\usepackage{amssymb}
\usepackage{amsmath}
\usepackage{amscd}
\usepackage{color}
\usepackage[T1]{fontenc}

\newcommand{\bg}{\begin{equation}}
\newcommand{\ed}{\end{equation}}
\newcommand{\bga}{\begin{eqnarray}}
\newcommand{\eda}{\end{eqnarray}}
\newcommand{\pf}{\textbf{Proof:\ }}

\def\cbdu{\par{\raggedleft$\Box$\par}}

\newtheorem {Theorem}  {Theorem}

\numberwithin{Theorem}{section}

\newtheorem {Lemma}[Theorem]  {Lemma}

\theoremstyle{definition}
\newtheorem{Definition}[Theorem]{Definition}
\theoremstyle{remark}

\expandafter\chardef\csname pre amssym.def
at\endcsname=\the\catcode`\@ \catcode`\@=11
\def\undefine#1{\let#1\undefined}
\def\newsymbol#1#2#3#4#5{\let\next@\relax
 \ifnum#2=\@ne\let\next@\msafam@\else
 \ifnum#2=\tw@\let\next@\msbfam@\fi\fi
 \mathchardef#1="#3\next@#4#5}
\def\mathhexbox@#1#2#3{\relax
 \ifmmode\mathpalette{}{\m@th\mathchar"#1#2#3}%
 \else\leavevmode\hbox{$\m@th\mathchar"#1#2#3$}\fi}
\def\hexnumber@#1{\ifcase#1 0\or 1\or 2\or 3\or 4\or 5\or 6\or 7\or 8\or
 9\or A\or B\or C\or D\or E\or F\fi}

\font\teneufm=eufm10 \font\seveneufm=eufm7 \font\fiveeufm=eufm5
\newfam\eufmfam
\textfont\eufmfam=\teneufm \scriptfont\eufmfam=\seveneufm
\scriptscriptfont\eufmfam=\fiveeufm

\catcode`\@=\csname pre amssym.def at\endcsname

\newcounter{remark}
\renewcommand{\a}{\alpha}

\newcommand{\s}{\sigma}

\newcommand{\R}{\mathbf{R}}

\renewcommand{\div}{\mbox{div}}

\def  \R   {{\mathbb R}}

\def  \12  {{\frac{1}{2}}}

\def\build#1_#2^#3{\mathrel{\mathop{\kern 0pt#1}\limits_{#2}^{#3}}}

\begin{document}

\title[Regularity of Chemotaxis-Navier-Stokes]{Low Modes Regularity criterion for a chemotaxis-Navier-Stokes system}


\author [Mimi Dai]{Mimi Dai}
\address{Department of Mathematics, Stat. and Comp. Sci.,  University of Illinois Chicago, Chicago, IL 60607,USA}
\email{mdai@uic.edu} 

\author [Han Liu]{Han Liu}
\address{Department of Mathematics, Stat. and Comp. Sci.,  University of Illinois Chicago, Chicago, IL 60607,USA}
\email{hliu94@uic.edu} 






\begin{abstract}

In this paper we study the regularity problem of a three dimensional chemotaxis-Navier-Stokes system on a periodic domain. A new regularity criterion in terms of only low modes of the oxygen concentration and the fluid velocity is obtained
 via a wavenumber splitting approach. The result improves many existing criteria in the literature.

\bigskip


\hspace{0.02cm}CLASSIFICATION CODE: 76D03, 35Q35, 35Q92, 92C17.
\end{abstract}

\maketitle

\section{Introduction}

We are interested in the following chemotaxis-Navier-Stokes system on a periodic domain
\begin{equation}\label{cmtns}
\begin{cases}
n_t+u\cdot \nabla n=\Delta n-\nabla \cdot(n\chi(c)\nabla c)\\
c_t+u \cdot \nabla c=\Delta c-nf(c)\\
u_t+(u\cdot \nabla)u+\nabla P =\Delta u +n \nabla \Phi\\
\nabla \cdot u=0, \ \ \ \ (t,x) \in \mathbb{R}^+ \times \mathbb{T}^3.
\end{cases}
\end{equation}
This coupled system arises from modelling aerobic bacteria, e.g. $\textit{Bacillus subtilis},$ suspended into sessile drops of water. It describes a scenario in which both the bacteria, whose population density is denoted by $n=n(t,x),$ and oxygen, whose concentration is denoted by $c=c(t,x),$ are transported by the fluid and diffuse randomly. In addition, the bacteria, which have chemotactic sensitivity  $\chi(c),$ tend to swim towards their nutrient oxygen and consume it at a per-capita rate $f(c)$.  At the same time, since the bacteria are heavier than water, their chemotactic swimming induces buoyant forces which affects the fluid motion. This buoyancy-driven effect is reflected in the third equation in system (\ref{cmtns}), represented by an extra term $n\nabla \Phi$ added to the Navier-Stokes equation. In this extra term, $\Phi$ denotes the gravitational potential, whereas the Navier-Stokes equation is conventionally written with $u=u(t,x)$ denoting the fluid velocity, and $P=P(t,x)$ the pressure. In this paper, we consider a simple yet prototypical case in which 
\begin{equation}\label{phchf}
\nabla\Phi \equiv \text{const.},\ \ \ \  \chi(c) \equiv \text{const.},\ \ \ \ f(c) \equiv c.
\end{equation} 
We note that in this case, solutions to system (\ref{cmtns}) satisfy the following scaling property: 
\[n_\lambda(t, x)=\lambda^2 n(\lambda^2t, \lambda x),\ \ c_\lambda(t, x)= c(\lambda^2t, \lambda x), \]
\[u_\lambda(t, x)=\lambda u(\lambda^2t, \lambda x),\ \ p_\lambda(t, x)=\lambda^2 p(\lambda^2t, \lambda x)\]
solve (\ref{cmtns}) with initial data 
$$n_{\lambda, 0}=\lambda^2n(\lambda x),\ \ c_{\lambda,0}= c(\lambda x),\ \ u_{\lambda, 0}=\lambda u(\lambda x),$$
provided that $$(n(t,x), c(t,x), u(t,x))$$ solves (\ref{cmtns}) with initial data $(n_0(x), c_0(x), u_0(x)).$ 
It is obvious that the Sobolev space $\dot H^{-\frac12}\times \dot H^{\frac12}\times\dot H^{\frac32}$ is scaling invariant (also called critical) for $(n,u,c)$ under the above nature scaling of the system.

Experiments showed that under the chemotaxis-fluid interaction of system (\ref{cmtns}), even almost homogeneous initial bacteria distribution can evolve and exhibit quite intricate spatial patterns, see \cite{DCCGK}, \cite{TCDWKG} and \cite{L1}. In \cite{L1} Lorz proved the existence of a local weak solution to the $3$D chemotaxis-Navier-Stokes system on bounded domains. In a recent work by Winkler (\cite{W2}), the existence of global weak solutions was proved under more general assumptions via entropy-energy estimates. We refer readers to the work of Winkler (\cite{W, W1, W2}), Liu and Lorz (\cite{LL}), Duan, Lorz and Markowich (\cite{DLM1}), Chae, Kang and Lee (\cite{CKL, CKL1}), as well as Jiang, Wu and Zheng (\cite{JWZ}) for more details about the well-posedness results for the chemotaxis-Navier-Stokes system.  

We are aware of several regularity criteria concerning the three dimensional chemotaxis-Navier-Stokes system. In \cite{CKL}, Chae, Kang and Lee obtained local-in-time classical solutions and Prodi-Serrin type regularity criteria. In particular, suppose that 
\begin{equation}\label{rc1}
\begin{split}
&\|u\|_{L^q(0,T; L^p)}+\|\nabla c\|_{L^2(0,T;L^\infty)}<\infty, \ \ \ \frac{3}{p}+\frac{2}{q}=1, \ \ 3 < p \leq \infty,
\end{split}
\end{equation} 
then the corresponding classical solution can be extended beyond time $T$. In \cite{CKL1}, Chae, Kang and Lee also obtained regularity criteria in terms of the $L^p$ norms of $u$ and $n.$ Jiang, Wu and Zheng, in their recent paper \cite{JWZ1}, proved that a classical solution to the initial boundary value problem of the Keller-Segel model i.e. the fluid free version of system (\ref{cmtns}), exists beyond time $T$ if 
\begin{equation}\notag
\begin{split}
&\|\nabla c\|_{L^2(0,T;L^\infty)}<\infty,\\
\text{or }\ \ &\|n\|_{L^q(0,T; L^p)}<\infty, \ \ \ \frac{3}{p}+\frac{2}{q} \leq 2,\ \ \frac{3}{2}<p \leq \infty.
\end{split}
\end{equation} 


In this paper, we aim to establish a regularity condition
weaker than (\ref{rc1}). Our result is stated in the following theorem, where $u_{\leq Q_u}$ and $c_{\leq Q_c}$ denote certain low modes of the velocity and oxygen concentration, which we shall explain later.  
\begin{Theorem}\label{Rgc}
Let $(n(t), c(t), u(t))$ be a weak solution to (\ref{cmtns}) on $[0,T]$ on the torus $\mathbb T^3$.
Assume that $(n(t), c(t), u(t))$ is regular on $[0, T)$ and
\begin{equation}\label{c1u1q}
\int^{T}_0 \|\nabla c_{\leq Q_c(t)}(t)\|^2_{L^ \infty}+\|u_{\leq Q_u(t)}(t)\|_{B^1_{\infty, \infty}} \mathrm{d}t < \infty,  
\end{equation}
then $(n(t), c(t), u(t))$ is regular on $[0, T].$
\end{Theorem}


We note that the quantity in (\ref{c1u1q}) is invariant with respect to the scaling of system (\ref{cmtns}).
It is obvious that the condition on the oxygen concentration $c$ in (\ref{c1u1q}) is weaker than that of (\ref{rc1}). It was also shown that the condition on velocity $u$ in (\ref{c1u1q}) is weaker than that of (\ref{rc1}), see \cite{CD1}. 

Our main devices are frequency localization technique and the wavenumber splitting method, 
which have been extensively applied to study the regularity problem and related problems (such as determining modes, see \cite{CD} and \cite{CDK}) for various supercritical dissipative equations.
In particular, in \cite{CS2} it was proved that a solution to the Navier-Stokes or Euler equation does not blow up at time $T$ if
\begin{gather*}
\int_0^T \|\nabla \times u_{\leq Q}\|_{B^0_{\infty, \infty}}\mathrm{d}t < \infty
\end{gather*}
for some wavenumber $\Lambda(t)=2^{Q(t)}.$ This result improved previous regularity criteria of Beale-Kato-Majda \cite{BKM} as well as of Prodi-Serrin-Ladyzhenskaya. The idea of separating high and low frequency modes by a critical wavenumber originats in Kolmogorov's theory of turbulence, which predicts the existence of a critical wavenumber above which the dissipation term dominates. 
Later, the wavenumber splitting method has been applied to a liquid crystal model with Q-tensor configuration \cite{D-Qtensor}, where it was shown that a condition solely on the low modes of the velocity (no condition on the Q-tensor) can guarantee regularity.


\bigskip

\section{Preliminaries}
\label{sec:pre}

\subsection{Notation}
\label{sec:notation}
The symbol $A\lesssim B$ denotes an estimate of the form $A\leq C B$ with
some absolute constant $C$, and $A\sim B$ denotes an estimate of the form $C_1
B\leq A\leq C_2 B$ with absolute constants $C_1$, $C_2$. 
The Sobolev norm $\|\cdot\|_{L^p}$ is shortened as $\|\cdot\|_p$ without no confusion. 
The symbols $W^{k,p}$ and $H^{s}$ represent the standard Sobolev spaces and $L^2$-based Sobolev spaces, respectively.

\subsection{Littlewood-Paley decomposition}
\label{sec:LPD}
The main analysis tools are the frequency localization method and a wavenumber splitting approach based on the Littlewood-Paley theory, which we briefly recall here. For a complete description of the theory and applications, the readers are referred to the books \cite{BCD} and \cite{G}.

We denote the Fourier transform and its inverse by $\mathcal F$ and $\mathcal F^{-1},$ respectively. We construct a family of smooth functions $\{\varphi_q \}_{q=-1}^\infty$ with annular support that forms a dyadic partition of unity in the frequency space, defined as 
\begin{equation}\notag
\varphi_q(\xi)=
\begin{cases}
\varphi(\lambda_q^{-1}\xi)  \ \ \ \mbox { for } q\geq 0,\\
\chi(\xi) \ \ \ \mbox { for } q=-1,
\end{cases}
\end{equation}
where $\varphi(\xi)=\chi(\xi/2)-\chi(\xi)$ and $\chi\in C_0^\infty(\R^n)$ is a nonnegative radial function chosen in a way such that 
\begin{equation}\notag
\chi(\xi)=
\begin{cases}
1, \ \ \mbox { for } |\xi|\leq\frac{3}{4}\\
0, \ \ \mbox { for } |\xi|\geq 1.
\end{cases}
\end{equation}

Introducing $\tilde h:=\mathcal F^{-1}\chi$ and $h:=\mathcal F^{-1}\varphi,$ we define the Littlewood-Paley projections for a function $u \in \mathcal{S}'$ as 
\begin{equation}\notag
\begin{cases}
& u_{-1}=\mathcal F^{-1}(\chi(\xi)\mathcal Fu)=\displaystyle\int \tilde h(y)u(x-y)dy,\\
&u_q:=\Delta_qu=\mathcal F^{-1}(\varphi(\lambda_q^{-1}\xi)\mathcal Fu)=\lambda_q^n\displaystyle\int h(\lambda_qy)u(x-y)dy, \ q\geq 0.
\end{cases}
\end{equation}
Then the identity
\bg\notag
u=\sum_{q=-1}^\infty u_q
\ed
holds in the sense of distributions. To simplify the notation, we denote
\bg\notag
\tilde u_q=u_{q-1}+u_q+u_{q+1}, \qquad u_{\leq Q}=\sum_{q=-1}^Qu_q,  \qquad u_{(P,Q]}=\sum_{q=P+1}^Qu_q.
\ed

We note that $ \|u\|_{H^s} \sim \left(\sum_{q=-1}^\infty\lambda_q^{2s}\|u_q\|_2^2\right)^{\frac{1}{2}},$
for each $u \in H^s$ and $s\in\R$. Utilizing the Littlewood-Paley theory, we give an equivalent definition of the norm of the Besov space $B_{p,\infty}^{s}$ as follows.
\begin{Definition}
Let $s\in \mathbb R$, and $1\leq p\leq \infty$. The Besov space $B_{p,\infty}^{s}$ is the space of tempered distributions $u$ whose Besov norm $\|u\|_{B_{p, \infty}^{s}} < \infty,$ where
$$\|u\|_{B_{p, \infty}^{s}}:=\sup_{q\geq -1}\lambda_q^s\|u_q\|_p.
$$
\end{Definition}

Moreover, we recall Bernstein's inequality, see \cite{L}.
\begin{Lemma}\label{brnst}
Let $n$ be the space dimension and $1\leq s\leq r\leq \infty$. Then for all tempered distributions $u$, 
\bg\notag
\|u_q\|_{r}\lesssim \lambda_q^{n(\frac{1}{s}-\frac{1}{r})}\|u_q\|_{s}.
\ed
\end{Lemma}
Throughout the paper, we will also utilize Bony's paraproduct decomposition
\begin{equation}\notag
\begin{split}
\Delta_q(u\cdot v)=&\sum_{|q-p|\leq 2}\Delta_q(u_{\leq{p-2}}\cdot v_p)+
\sum_{|q-p|\leq 2}\Delta_q(u_{p}\cdot v_{\leq{p-2}})+\sum_{p\geq q-2}\Delta_q(u_p\cdot\tilde v_p),
\end{split}
\end{equation}
as well as the commutator notation 
$$
[\Delta_q, u_{\leq{p-2}}\cdot\nabla]v_p=\Delta_q(u_{\leq{p-2}}\cdot\nabla v_p)-u_{\leq{p-2}}\cdot\nabla \Delta_qv_p.
$$
An estimate for the commutator is given by the following lemma (c.f. \cite{CD}).
\begin{Lemma}\label{comtr}
Let $\frac1{r_2}+\frac1{r_3}=\frac1{r_1},$ we have the estimate
\begin{equation}\notag
\|[\Delta_q,u_{\leq{p-2}}\cdot\nabla] v_q\|_{r_1}\lesssim \|v_q\|_{r_3}\sum_{p' \leq p-2} \lambda_{p'} \|u_{p'}\|_{r_2}.
\end{equation}
\end{Lemma}

\subsection{Weak solution and regular solution to system (\ref{cmtns})}
\label{sec:sol} 
From \cite{W2}, we know that on bounded, smooth and convex domain $\Omega$ in three dimension, system (\ref{cmtns}) has a global weak solution $(n,c,u)$ which satisfies the equations in (\ref{cmtns}) in the distributional sense, provided that the initial data $(n_0, c_0, u_0)$ satisfy
$$n_0 \in L \log L(\Omega), \ \ c_0 \in L^\infty(\Omega) \text{ and } \sqrt{c_0} \in H^1(\Omega), \ \ u_0 \in L_\s^2(\Omega),$$ where $L_\s^2(\Omega)$ is the space of solenoidal vector field in $L^2(\Omega).$  We highlight the following properties of the weak solution $(n,c,u)$ in particular
\begin{align*}
&n \in L^\infty(0,\infty; L^1(\Omega)),\ \ c \in L^\infty(0,\infty; L^\infty(\Omega)),\\
&u \in L^\infty_{loc}(0,\infty; L^2(\Omega))\cap L^2_{loc}(0,\infty; H^1_0(\Omega)).
\end{align*}
A regular solution of (\ref{cmtns}) is understood in the way that the solution has enough regularity to satisfy the equations of the system point-wise. Typically, a solution in a space with higher regularity than its critical space can be shown regular via bootstrapping arguments.  The local existence of regular solution to (\ref{cmtns}) was shown in \cite{CKL}.


\subsection{Parabolic regularity theory}
\label{sec:par}
We consider the heat equation on $\mathbb{T}^d$   with $d\geq2$
\begin{equation}\label{eq-heat}
u_t-\Delta u=f
\end{equation}
with initial data $u_0$.
We shall see that the solution $u$ turns out to be smoother than the source term $f.$
\begin{Lemma}\label{le-heat}
Let $u$ be a solution to (\ref{eq-heat}) with $u_0\in H^{\alpha+1}$ and $f\in L^2(0,T; H^\alpha)$ for $\alpha\in \mathbb R$. Then we have $u\in L^2(0,T; H^{\alpha+2})\cap H^1(0,T; H^\alpha)$.
\end{Lemma}
\pf
Projecting equation (\ref{eq-heat}) by $\Delta_q$ and taking inner product of the resulted equation with $\lambda_q^{2\alpha+4}u_q$ leads 
\begin{equation}\notag
\frac12\frac{d}{dt}\lambda_q^{2\alpha+4}\|u_q\|_2^2+\lambda_q^{2\alpha+4}\|\nabla u_q\|_2^2=\lambda_q^{2\alpha+4}\int f_qu_q\, \mathrm{d}x.
\end{equation}
Applying H\"older's and Young's inequalities to the right hand side yields
\begin{equation}\notag
\frac{d}{dt}\lambda_q^{2\alpha+4}\|u_q\|_2^2+\lambda_q^{2\alpha+4}\|\nabla u_q\|_2^2\leq 4\lambda_q^{2\alpha+2}\|f_q\|_2^2.
\end{equation}
As a consequence of Duhamel's formula, summation in $q$ and integration over $[0,T]$, we obtain
\begin{equation}\notag
\begin{split}
\int_0^T\sum_{q \geq -1}\lambda_q^{2\a+4}\|u_q(t)\|_2^2\, \mathrm{d}t \leq &\int_0^T\sum_{q \geq -1} \lambda_q^{2\a+4}\|u_q(0)\|_2^2e^{-\lambda_q^2t}\, \mathrm{d}t\\
&+4\int_0^T\sum_{q \geq -1}\lambda_q^{2\a+2}\int_0^te^{-\lambda_q^2(t-s)}\|f_q(s)\|_2^2\mathrm{d}s \, \mathrm{d}t.\
\end{split}
\end{equation}
The first integral on the right hand side is handled as
\begin{equation}\notag
\int_0^T\sum_{q \geq -1} \lambda_q^{2\a+4}\|u_q(0)\|_2^2e^{-\lambda_q^2t}\, \mathrm{d}t
\leq  \sum_{q \geq -1} \lambda_q^{2\a+2}\|u_q(0)\|_2^2\left(1-e^{-\lambda_q^2T}\right)
\lesssim \|u_0\|_{H^{\alpha+1}}^2.
\end{equation}
In order to estimate the second integral, we exchange the order of integration to obtain
\begin{equation}\notag
\begin{split}
&\int_0^T\sum_{q \geq -1}\lambda_q^{2\a+2}\int_0^te^{-\lambda_q^2(t-s)}\|f_q(s)\|_2^2\mathrm{d}s \, \mathrm{d}t\\
\leq &\int_0^T\int_s^T\sum_{q \geq -1}\lambda_q^{2\a+2}e^{-\lambda_q^2(t-s)}\|f_q(s)\|_2^2\mathrm{d}t \, \mathrm{d}s\\
\leq &\int_0^T\sum_{q \geq -1}\lambda_q^{2\a}\|f_q(s)\|_2^2\left(1-e^{-\lambda_q^2(T-s)}\right)\, \mathrm{d}s\\
\lesssim &\|f\|_{L^2(0,T;H^\alpha)}^2.
\end{split}
\end{equation}
Combining the estimates above, we conclude that $u\in L^2(0,T; H^{\alpha+2})$ for  $\alpha\in \mathbb R$.

To prove $u\in H^1(0,T; H^{\alpha})$, we first project equation (\ref{eq-heat}) to the $q$-th dyadic shell
\[(u_t)_q=\Delta u_q+f_q.\]
It follows 
\[\|(u_t)_q\|_2^2\leq 2\|\Delta u_q\|_2^2+2\|f_q\|_2^2.\]
Thus we deduce that
\begin{equation}\notag
\begin{split}
\int_0^T\sum_{q\geq -1}\lambda_q^{2\alpha} \|(u_t)_q\|_2^2\,\mathrm{d}t\lesssim &
\int_0^T\sum_{q\geq -1}\lambda_q^{2\alpha} \|\Delta u_q\|_2^2\,\mathrm{d}t
+\int_0^T\sum_{q\geq -1}\lambda_q^{2\alpha} \|f_q\|_2^2\,\mathrm{d}t\\
\lesssim &\|u\|_{L^2(0,T; H^{\alpha+2})}^2+\|f\|_{L^2(0,T;H^\alpha)}^2.
\end{split}
\end{equation}
It is then clear that $u\in H^1(0,T; H^\alpha)$, which completes the proof of the lemma.

\cbdu


\bigskip

\section{Proof of Theorem \ref{Rgc}}
\label{sec:reg}

This section is devoted to the proof of the main result.
We start by introducing the dissipation wavenumber $\Lambda_u(t)$ for $u$ and $\Lambda_c(t)$ for $c$, 
\begin{equation}\label{wave-uc}
\begin{split}
\Lambda_u(t)=&\min\left\{\lambda_q: \lambda_p^{-1}\|u_p(t)\|_\infty<C_0, \forall p >q, q \in \mathbb{N}\right\},\\
\Lambda_c(t)=&\min\left\{\lambda_q: \lambda_p^{\frac{3}{r}}\|c_p(t)\|_r<C_0, \forall p >q, q \in \mathbb{N}\right\}, \ \ r\in\left(3,\frac{3}{1-\varepsilon}\right), 
\end{split}
\end{equation}
where $\varepsilon>0$ is a fixed arbitrarily small constant, and $C_0$ is a small constant to be determined later. Through out this section, we use $C$ for various absolute constants which can be different from line to line. 
In addition, we let $Q_u(t)$ and $Q_c(t)$ be integers such that $\lambda_{Q_u(t)}=\Lambda_u(t)$ and $\lambda_{Q_c(t)}=\Lambda_c(t)$. Then the constraint on the low modes is defined as

\[f(t):= \| \nabla c_{\leq Q_c(t)}(t)\|^2_{L^{\infty}}+\|u_{\leq Q_u(t)}(t)\|_{B^1_{\infty, \infty}}.\]
Notice that the wavenumber $\Lambda_u$ separates the inertial range from the dissipation range where the viscous term $\Delta u$ dominates; and $\Lambda_c$ has the same property.
Precisely, we have
\[\|u(t)_{Q_u}\|_{\infty}\geq C_0 \Lambda_u(t), \ \   \Lambda_c^{\frac3r}(t)\|c(t)_{Q_c}\|_{r}\geq C_0;\]
\[\lambda_q\|u(t)_q\|_{\infty}< C_0, \ \forall q>Q_u; \ \   \lambda_q^{\frac3r}(t)\|c(t)_q\|_{r}< C_0, \ \forall q>Q_c.\]

The crucial part is to establish a uniform (in time) bound for each of the unknowns $n, u$ and $c$ in a space with higher regularity than the critical Sobolev space. In fact, it is sufficient to prove that $(n, u, c)\in L^\infty(0,T; \dot H^{s_1})\times L^\infty(0,T; \dot H^{s_2})\times L^\infty(0,T; \dot H^{s_3})$ for some $s_1>-\frac12, s_2>\frac12$ and $s_3>\frac32$. Due to the complicated interactions among the three equations in (\ref{cmtns}), the aforementioned goal will be achieved in two steps. The first step is to show that $(n, u, c)\in L^\infty(0,T; \dot H^{s})\times L^\infty(0,T; \dot H^{s+1})\times L^\infty(0,T; \dot H^{s+1})$ for some $s\in (-\frac12,0)$.
The second step is to apply bootstrapping arguments, the Lp-Lq theory for parabolic equations and a mixed derivative theorem to the equation of oxygen concentration $c$, and hence improve the regularity of $c$.

In the first step, we multiply the equations in (\ref{cmtns}) by $\lambda_q^{2s}\Delta_q^2n,$ $\lambda_q^{2s+2}\Delta_q^2c$ and $\lambda_q^{2s+2}\Delta_q^2u,$ respectively. Integrating and summing lead to 

\begin{equation}\label{eqtnn}
\begin{split}
\frac{1}{2}\frac{\mathrm{d}}{\mathrm{d}t}\sum_{q \geq -1} \lambda^{2s}_q \|n_q\|_2^2 \leq & -\sum_{q \geq -1}\lambda^{2s}_q\|\nabla n_q\|_2^2 -\sum_{q \geq -1}\lambda^{2s}_q\int_{\mathbb{R}^3} \Delta_q (u \cdot \nabla n)n_q \mathrm{d}x\\
&-\sum_{q \geq -1}\lambda^{2s}_q\int_{\mathbb{R}^3} \Delta_q(\nabla \cdot (n \chi(c)\nabla c))n_q \mathrm{d}x;
\end{split}
\end{equation}
\begin{equation}\label{eqtnc}
\begin{split}
&\frac{1}{2}\frac{\mathrm{d}}{\mathrm{d}t}\sum_{q \geq -1} \lambda^{2s+2}_q \|c_q\|_2^2 \leq  -\sum_{q \geq -1}\lambda^{2s+2}_q\|\nabla c_q\|_2^2 \\
&-\sum_{q \geq -1}\lambda^{2s+2}_q\int_{\mathbb{R}^3} \Delta_q (u \cdot \nabla c)c_q \mathrm{d}x
-\sum_{q \geq -1}\lambda^{2s+2}_q\int_{\mathbb{R}^3} \Delta_q(nf(c))c_q \mathrm{d}x; 
\end{split}
\end{equation}
\begin{equation}\label{eqtnu}
\begin{split}
&\frac{1}{2}\frac{\mathrm{d}}{\mathrm{d}t}\sum_{q \geq -1} \lambda^{2s+2}_q \|u_q\|_2^2 \leq  -\sum_{q \geq -1}\lambda^{2s+2}_q\|\nabla u_q\|_2^2 \\
&-\sum_{q \geq -1}\lambda^{2s+2}_q\int_{\mathbb{R}^3} \Delta_q (u \cdot \nabla u)u_q \mathrm{d}x
-\sum_{q \geq -1}\lambda^{2s+2}_q\int_{\mathbb{R}^3} \Delta_q(n\nabla\Phi)u_q \mathrm{d}x.
\end{split}
\end{equation}
For simplicity we label the terms
\begin{gather*}
I:=-\sum_{q \geq -1}\lambda^{2s+2}_q\int_{\mathbb{R}^3}\Delta_q (u \cdot \nabla u)u_q \mathrm{d}x, \ \ \ \ II:=-\sum_{q \geq -1}\lambda^{2s+2}_q\int_{\mathbb{R}^3} \Delta_q (u \cdot \nabla c)c_q \mathrm{d}x,\\
III:= -\sum_{q \geq -1}\lambda^{2s}_q\int_{\mathbb{R}^3} \Delta_q (u \cdot \nabla n)n_q \mathrm{d}x,\\
IV:=-\sum_{q \geq -1}\lambda^{2s+2}_q\int_{\mathbb{R}^3} \Delta_q(n\nabla\Phi)u_q \mathrm{d}x, \ \ \ \ V:=-\sum_{q \geq -1}\lambda^{2s+2}_q\int_{\mathbb{R}^3} \Delta_q(nf(c))c_q \mathrm{d}x,\\
VI:=-\sum_{q \geq -1}\lambda^{2s}_q\int_{\mathbb{R}^3} \Delta_q(\nabla \cdot (n \chi(c)\nabla c))n_q \mathrm{d}x.
\end{gather*}

\subsection{Estimate of $I$}
We estimate the term $I$ using the wavenumber splitting method. As we shall see, the commutator reveals certain cancellation within the nonlinear interactions. 
\begin{Lemma}\label{ppsuu}
Let $s>-\frac{1}{2}.$ We have
\begin{equation}\notag
|I| \lesssim C_0\sum_{q> -1}\lambda_q^{2s+4}\|u_q\|_2^2+Q_uf(t)\sum_{q\geq -1}\lambda_q^{2s+2}\|u_q\|_2^2.
\end{equation} 
\end{Lemma}
\pf
Applying Bony's paraproduct decomposition to $I$ leads to
\begin{equation}\notag
\begin{split}
I=&-\sum_{q\geq -1}\sum_{|q-p|\leq 2}\lambda_q^{2s+2}\int_{\R^3}\Delta_q(u_{\leq p-2}\cdot\nabla u_{p})u_q\, \mathrm {d}x\\
&-\sum_{q\geq -1}\sum_{|q-p|\leq 2}\lambda_q^{2s+2}\int_{\R^3}\Delta_q(u_{p}\cdot\nabla u_{\leq{p-2}})u_q\, \mathrm {d}x\\
&-\sum_{q\geq -1}\sum_{p\geq q-2}\lambda_q^{2s+2}\int_{\R^3}\Delta_q(u_p\cdot\nabla\tilde u_p)u_q\, \mathrm {d}x\\
=:&I_{1}+I_{2}+I_{3}.
\end{split}
\end{equation}

Using the fact $\sum_{|q-p|\leq 2}\Delta_q u_p=u_q$ and the commutator notation, we have
\begin{equation}\notag
\begin{split}
I_{1}=&-\sum_{q\geq -1}\sum_{|q-p|\leq 2}\lambda_q^{2s+2}\int_{\R^3}[\Delta_q, u_{\leq{p-2}}\cdot\nabla] u_pu_q\, \mathrm {d}x\\
&-\sum_{q\geq -1}\lambda_q^{2s_3}\int_{\R^3}u_{\leq{q-2}}\cdot\nabla  u_q u_q\, \mathrm {d}x\\
&-\sum_{q\geq -1}\sum_{|q-p|\leq 2}\lambda_q^{2s+2}\int_{\R^3}(u_{\leq{p-2}}-u_{\leq{q-2}})\cdot\nabla\Delta_qu_p u_q\, \mathrm {d}x\\
=:&I_{11}+I_{12}+I_{13}.
\end{split}
\end{equation}
Moreover we have $I_{12}=0$ due to that $\div\, u_{\leq q-2}=0.$ 

We then split $I_{11}$ based on definition of $\Lambda_u(t)$
\begin{equation}\notag
\begin{split}
|I_{11}|\leq &\sum_{q\geq -1}\sum_{|q-p|\leq 2}\lambda_q^{2s+2}\int_{\R^3}\left|[\Delta_q, u_{\leq{p-2}}\cdot\nabla] u_pu_q\right|\, \mathrm {d}x\\
\leq & \sum_{p\leq Q_u+2}\sum_{|q-p|\leq 2}\lambda_q^{2s+2}\int_{\R^3}\left|[\Delta_q, u_{\leq{p-2}}\cdot\nabla] u_pu_q\right|\, \mathrm {d}x\\
&+\sum_{p> Q_u+2}\sum_{|q-p|\leq 2}\lambda_q^{2s+2}\int_{\R^3}\left|[\Delta_q, u_{\leq Q_u}\cdot\nabla] u_pu_q\right|\, \mathrm {d}x\\
&+\sum_{p> Q_u+2}\sum_{|q-p|\leq 2}\lambda_q^{2s+2}\int_{\R^3}\left|[\Delta_q, u_{(Q_u,p-2]}\cdot\nabla] u_pu_q\right|\, \mathrm {d}x\\
=: &I_{111}+I_{112}+I_{113}.
\end{split}
\end{equation}
Using (\ref{comtr}), H\"older's inequality, and  definition of $f(t)$, we obtain
\begin{equation}\notag
\begin{split}
I_{111}\leq &\sum_{1\leq p\leq Q_u+2}\sum_{|q-p|\leq 2}\lambda_q^{2s+2}\|\nabla u_{\leq p-2}\|_\infty\|u_p\|_2\|u_q\|_2\\
\lesssim & f(t) \sum_{1\leq p\leq Q_u+2}\|u_p\|_2\sum_{|q-p|\leq 2}\lambda_q^{2s+2}\|u_q\|_2\sum_{p'\leq p-2}1\\
\lesssim &Q_u f(t) \sum_{1\leq p\leq Q_u+2}\lambda_p^{s+1}\|u_p\|_2\sum_{|q-p|\leq 2}\lambda_q^{s+1}\|u_q\|_2\\
\lesssim &Q_u f(t) \sum_{q\geq -1}\lambda_q^{2s+2}\|u_q\|_2^2;
\end{split} 
\end{equation}
and similarly
\begin{equation}\notag
\begin{split}
I_{112}\leq &\sum_{ p> Q_u+2}\sum_{|q-p|\leq 2}\lambda_q^{2s+2}\|\nabla u_{\leq Q_u}\|_\infty\|u_p\|_2\|u_q\|_2\\
\lesssim & Q_uf(t) \sum_{ p> Q_u+2}\|u_p\|_2\sum_{|q-p|\leq 2}\lambda_q^{2s+2}\|u_q\|_2\\
\lesssim & Q_uf(t) \sum_{ p> Q_u+2}\lambda_p^{s+1}\|u_p\|_2\sum_{|q-p|\leq 2}\lambda_q^{s+1}\|u_q\|_2\\
\lesssim &Q_u f(t) \sum_{q> Q_u}\lambda_q^{2s+2}\|u_q\|_2^2.
\end{split}
\end{equation}
We estimate $I_{113}$ with the help of H\"older's inequality and Lemma (\ref{comtr}) 
\begin{equation}\notag
\begin{split}
I_{113}
\leq &\sum_{ p> Q_u+2}\sum_{|q-p|\leq 2}\lambda_q^{2s+2}\|[\Delta_q,u_{(Q_u,p-2]}\cdot\nabla]u_p\|_2\|u_q\|_2\\
\leq &\sum_{ p> Q_u+2}\|u_p\|_2\sum_{|q-p|\leq 2}\lambda_q^{2s+2}\|u_q\|_2\sum_{Q_u<p'\leq p-2}\lambda_{p'}\|u_{p'}\|_\infty\\
\lesssim & C_0 \sum_{ p> Q_u+2}\|u_p\|_2\sum_{|q-p|\leq 2}\lambda_q^{2s+2}\|u_q\|_2\sum_{Q_u<p'\leq p-2}\lambda_{p'}^{2}\\
\lesssim & C_0 \sum_{ p> Q_u+2}\lambda_p^{2s+4}\|u_p\|_2^2\sum_{Q_u<p'\leq p-2}\lambda_{p'-p}^{2}\\
\lesssim & C_0 \sum_{ p> Q_u+2}\lambda_p^{2s+4}\|u_p\|_2^2.
\end{split}
\end{equation}

We split $I_{13}$ according to the definition of $\Lambda_u(t)$
\begin{equation}\notag
\begin{split}
|I_{13}|\leq & \sum_{q\geq -1}\sum_{|q-p|\leq 2}\lambda_q^{2s+2}\int_{\R^3}\big|(u_{\leq{p-2}}-u_{\leq{q-2}})\cdot\nabla\Delta_qu_p u_q\big|\, \mathrm {d}x\\
\leq & \sum_{-1\leq q\leq Q_u}\sum_{|q-p|\leq 2}\lambda_q^{2s+2}\int_{\R^3}\big|(u_{\leq{p-2}}-u_{\leq{q-2}})\cdot\nabla\Delta_qu_p u_q\big|\, \mathrm {d}x\\
 & +\sum_{q>Q_u}\sum_{|q-p|\leq 2}\lambda_q^{2s+2}\int_{\R^3}\big|(u_{\leq{p-2}}-u_{\leq{q-2}})\cdot\nabla\Delta_qu_p u_q\big|\, \mathrm {d}x\\
=:& I_{131}+I_{132}.
\end{split}
\end{equation}
Using H\"older's inequality and definition of $f(t)$ we can bound $I_{131}.$
\begin{equation}\notag
\begin{split}
|I_{131}|\leq & \sum_{-1\leq q\leq Q_u}\lambda_q^{2s+3} \|u_q\|_{\infty}\sum_{|q-p|\leq 2}\|u_{\leq{p-2}}-u_{\leq{q-2}}\|_2\|u_p\|_2\\
\lesssim & f(t)\sum_{-1\leq q\leq Q_u}\lambda_q^{2s+2} \sum_{|q-p|\leq 2}\|u_{\leq{p-2}}-u_{\leq{q-2}}\|_2\|u_p\|_2\\
\lesssim & f(t)\sum_{-1\leq q\leq Q_u}\lambda_q^{2s+2} \|u_q\|_2^2.\\
\end{split}
\end{equation}
And we estimate $I_{132}$ using H\"older's inequality and the definition of $\Lambda_u(t)$,
\begin{equation}\notag
\begin{split}
|I_{132}|\leq  & \sum_{q> Q_u}\lambda_q^{2s+3} \|u_q\|_\infty\sum_{|q-p|\leq 2}\|u_{\leq{p-2}}-u_{\leq{q-2}}\|_2\|u_p\|_{2}\\
\lesssim & C_0 \sum_{q> Q_u}\lambda_q^{2s+4}\sum_{|q-p|\leq 2}\|u_{\leq{p-2}}-u_{\leq{q-2}}\|_2\|u_p\|_2\\
\lesssim & C_0 \sum_{q> Q_u}\lambda_q^{2s+4}\|u_q\|_2^2.  
\end{split}
\end{equation}

We omit the detailed estimation of $I_{2}$ as it is similar to that of $I_{11}$. Meanwhile, for $I_{3}$ we have for $s>-\frac12$
\begin{equation}\notag
\begin{split}
|I_{3}|\leq&\sum_{q\geq -1}\sum_{p\geq q-2}\lambda_q^{2s+2}\int_{\mathbb R^3}|\Delta_q(u_p\otimes\tilde u_p)\nabla u_q|\, \mathrm {d}x\\
\leq &\sum_{q> Q_u}\lambda_q^{2s+3}\|u_q\|_\infty\sum_{p\geq q-2}\|u_p\|_2^2
+\sum_{-1\leq q\leq Q_u}\lambda_q^{2s+3}\|u_q\|_\infty\sum_{p\geq q-2}\|u_p\|_2^2\\
\leq &C_0\sum_{q> Q_u}\lambda_q^{2s+4}\sum_{p\geq q-2}\|u_p\|_2^2
+f(t)\sum_{-1\leq q\leq Q_u}\lambda_q^{2s+2}\sum_{p\geq q-2}\|u_p\|_2^2\\
\lesssim & C_0 \sum_{p> Q_u}\lambda_p^{2s+4}\|u_p\|_2^2\sum_{Q_u< q\leq p+2}\lambda_{q-p}^{2s+4}
+f(t)\sum_{p\geq -1}\lambda_p^{2s+2}\|u_p\|_2^2\sum_{q\leq p+2}\lambda_{q-p}^{2s+2}\\
\lesssim & C_0 \sum_{q> Q_u}\lambda_q^{2s+4}\|u_q\|_2^2+f(t)\sum_{q\geq -1}\lambda_q^{2s+2}\|u_q\|_2^2.
\end{split}
\end{equation}
We combine the above estimates to conclude that
\begin{equation}\notag
|I|\lesssim C_0\sum_{q> -1}\lambda_q^{2s+4}\|u_q\|_2^2+Q_uf(t)\sum_{q\geq -1}\lambda_q^{2s+2}\|u_q\|_2^2.
\end{equation}
\cbdu

\subsection{Estimate of $II$} 
\begin{Lemma}\label{ppscv}
Let $s>-\frac12$. We have
\begin{equation}\notag
\begin{split}
|II|\leq & \left(CC_0+\frac1{32}\right) \sum_{q>-1}\left(\lambda_q^{2s+4}\|c_q\|_2^2+\lambda_q^{2s+4}\|u_q\|_2^2\right)+CQ_uf(t)\sum_{q\geq-1}\lambda_q^{2s+2}\|c_q\|_2^2.
\end{split}
\end{equation}
\end{Lemma}
\pf
We start with Bony's paraproduct decomposition for $II$,
\begin{equation}\notag
\begin{split}
II=&\sum_{q\geq -1}\sum_{|q-p|\leq 2}\lambda_q^{2s+2}\int_{\R^3}\Delta_q(u_{\leq p-2}\cdot\nabla c_p)c_q\, \mathrm{d}x\\
&+\sum_{q\geq -1}\sum_{|q-p|\leq 2}\lambda_q^{2s+2}\int_{\R^3}\Delta_q(u_{p}\cdot\nabla c_{\leq{p-2}})c_q\, \mathrm{d}x\\
&+\sum_{q\geq -1}\sum_{p\geq q-2}\lambda_q^{2s+2}\int_{\R^3}\Delta_q(u_p\cdot\nabla\tilde c_p)c_q\, \mathrm{d}x\\
=&II_{1}+II_{2}+II_{3}.
\end{split}
\end{equation}
Using the commutator notation, we rewrite $II_1$ as 
\begin{equation}\notag
\begin{split}
II_{1}=&\sum_{q\geq -1}\sum_{|q-p|\leq 2}\lambda_q^{2s+2}\int_{\R^3}[\Delta_q,u_{\leq{p-2}}\cdot\nabla] c_pc_q\, \mathrm{d}x\\
&+\sum_{q\geq -1}\lambda_q^{2s+2}\int_{\R^3}u_{\leq{q-2}}\cdot\nabla c_q c_q\, \mathrm{d}x\\
&+\sum_{q\geq -1}\sum_{|q-p|\leq 2}\lambda_q^{2s+2}\int_{\R^3}(u_{\leq{p-2}}-u_{\leq{q-2}})\cdot\nabla\Delta_qc_p c_q\, \mathrm{d}x\\
=:&II_{11}+II_{12}+II_{13}.
\end{split}
\end{equation}
Here, just as in proposition (\ref{ppsuu}), we used $\sum_{|q-p|\leq 2}\Delta_pc_q=c_q$ to obtain $II_{12},$ which vanishes as $\div \, u_{\leq q-2}=0$. 

One can see that $II_{11}$ enjoys the same estimates that $I_{11}$ does. Splitting the summation by $Q_u(t)$ yields 
\begin{equation}\notag
\begin{split}
|II_{11}|\leq &\sum_{q\geq -1}\sum_{|q-p|\leq 2}\lambda_q^{2s+2}\int_{\R^3}\left|[\Delta_q, u_{\leq{p-2}}\cdot\nabla] c_pc_q\right|\, \mathrm {d}x\\
\leq & \sum_{p\leq Q_u+2}\sum_{|q-p|\leq 2}\lambda_q^{2s+2}\int_{\R^3}\left|[\Delta_q, u_{\leq{p-2}}\cdot\nabla] c_pc_q\right|\, \mathrm {d}x\\
&+\sum_{p> Q_u+2}\sum_{|q-p|\leq 2}\lambda_q^{2s+2}\int_{\R^3}\left|[\Delta_q, u_{\leq Q_u}\cdot\nabla] c_pc_q\right|\, \mathrm {d}x\\
&+\sum_{p> Q_u+2}\sum_{|q-p|\leq 2}\lambda_q^{2s+2}\int_{\R^3}\left|[\Delta_q, u_{(Q_u,p-2]}\cdot\nabla] c_pc_q\right|\, \mathrm {d}x\\
=: &II_{111}+II_{112}+II_{113}.
\end{split}
\end{equation}
We estimate the first two of the above three terms by H\"older's inequality, Lemma (\ref{comtr}) and the definition of $f(t).$
\begin{equation}\notag
\begin{split}
II_{111}\leq &\sum_{1\leq p\leq Q_u+2}\sum_{|q-p|\leq 2}\lambda_q^{2s+2}\|\nabla u_{\leq p-2}\|_\infty\|c_p\|_2\|c_q\|_2\\
\lesssim & Q_uf(t) \sum_{1\leq p\leq Q_u+2}\lambda_p^{s+1}\|c_p\|_2\sum_{|q-p|\leq 2}\lambda_q^{s+1}\|c_q\|_2\\
\lesssim &Q_u f(t) \sum_{q\geq -1}\lambda_q^{2s+2}\|c_q\|_2^2,
\end{split} 
\end{equation}
\begin{equation}\notag
\begin{split}
II_{112}\leq &\sum_{ p> Q_u+2}\sum_{|q-p|\leq 2}\lambda_q^{2s+2}\|\nabla u_{\leq Q_u}\|_\infty\|c_p\|_2\|c_q\|_2\\
\lesssim &Q_u f(t) \sum_{ p> Q_u+2}\lambda_p^{s+1}\|c_p\|_2\sum_{|q-p|\leq 2}\lambda_q^{s+1}\|c_q\|_2\\
\lesssim &Q_u f(t) \sum_{q> Q_u}\lambda_q^{2s+2}\|c_q\|_2^2.
\end{split}
\end{equation}
And with the help of H\"older's inequality, Lemma (\ref{comtr}) and the definition of $\Lambda_u(t)$, we estimate $I_{113}$ as
\begin{equation}\notag
\begin{split}
II_{113}
\leq &\sum_{ p> Q_u+2}\sum_{|q-p|\leq 2}\lambda_q^{2s+2}\|[\Delta_q,u_{(Q_u,p-2]}\cdot\nabla]c_p\|_2\|c_q\|_2\\
\leq &\sum_{ p> Q_u+2}\sum_{|q-p|\leq 2}\lambda_q^{2s+2}\|c_q\|_2\sum_{Q_u<p'\leq p-2}\lambda_{p'}\|u_{p'}\|_\infty\|c_p\|_2\\
\leq& \sum_{ p> Q_u+2}\|c_p\|_2\sum_{|q-p|\leq 2}\lambda_q^{2s+2}\|c_q\|_2\sum_{Q_u<p'\leq p-2}\lambda_{p'}\|u_{p'}\|_\infty\\
\lesssim & C_0 \sum_{ p> Q_u+2}\|c_p\|_2\sum_{|q-p|\leq 2}\lambda_q^{2s+2}\|c_q\|_2\sum_{Q_u<p'\leq p-2}\lambda_{p'}^{2}\\
\lesssim & C_0 \sum_{ p> Q_u+2}\lambda_p^{2s+4}\|c_p\|_2^2\sum_{Q_u<p'\leq p-2}\lambda_{p'-p}^{2}\\
\lesssim & C_0 \sum_{ p> Q_u+2}\lambda_p^{2s+4}\|c_p\|_2^2.
\end{split}
\end{equation}

$II_{13}$ can be estimated in the same fashion as $I_{13}.$ Splitting the sum by the wavenumber $Q_u(t),$ we have  
\begin{equation}\notag
\begin{split}
|II_{13}|\leq & \sum_{q\geq -1}\sum_{|q-p|\leq 2}\lambda_q^{2s+2}\int_{\R^3}\big|(u_{\leq{p-2}}-u_{\leq{q-2}})\cdot\nabla\Delta_qc_p c_q\big|\, \mathrm {d}x\\
\leq & \sum_{-1\leq q\leq Q_u}\sum_{|q-p|\leq 2}\lambda_q^{2s+2}\int_{\R^3}\big|(u_{\leq{p-2}}-u_{\leq{q-2}})\cdot\nabla\Delta_qc_p c_q\big|\, \mathrm {d}x\\
 & +\sum_{q>Q_u}\sum_{|q-p|\leq 2}\lambda_q^{2s+2}\int_{\R^3}\big|(u_{\leq{p-2}}-u_{\leq{q-2}})\cdot\nabla\Delta_qc_p c_q\big|\, \mathrm {d}x\\
=:&II_{131}+II_{132},
\end{split}
\end{equation}
which, using H\"older's inequality along with the  definition of $f(t)$ and $\Lambda_u(t),$ we estimate as
\begin{equation}\notag
\begin{split}
|II_{131}|\leq & \sum_{-1\leq q\leq Q_u}\lambda_q^{2s+3} \|c_q\|_2\sum_{|q-p|\leq 2}\|u_{\leq{p-2}}-u_{\leq{q-2}}\|_\infty\|c_p\|_2\\
\lesssim & f(t)\sum_{-1\leq q\leq Q_u}\lambda_q^{2s+2}\|c_q\|_2 \sum_{|q-p|\leq 2}\|c_p\|_2\\
\lesssim & f(t)\sum_{-1\leq q\leq Q_u}\lambda_q^{2s+2} \|c_q\|_2^2,\\
\end{split}
\end{equation}
\begin{equation}\notag
\begin{split}
|II_{132}|\leq  & \sum_{q> Q_u}\lambda_q^{2s+3} \|c_q\|_2\sum_{|q-p|\leq 2}\|u_{\leq{p-2}}-u_{\leq{q-2}}\|_\infty\|c_p\|_2\\
\lesssim & C_0 \sum_{q> Q_u}\lambda_q^{2s+3}\|c_q\|_2\sum_{|q-p|\leq 2}\lambda_p\|c_p\|_2\\
\lesssim & C_0 \sum_{q> Q_u}\lambda_q^{2s+4}\|c_q\|_2^2.  
\end{split}
\end{equation}

To estimate $II_2$, we make use of the wavenumber $Q_c(t)$ instead of $Q_u(t)$. Splitting the summation by $Q_c(t)$ yields
\begin{equation}\notag
\begin{split}
II_2=&  \sum_{q \geq -1}\sum_{|q-p| \leq 2}\lambda^{2s+2}_q \int_{\mathbb{R}^3}\Delta_q (u_p\cdot \nabla c_{\leq p-2})c_q\mathrm{d}x\\
=&  \sum_{q \geq -1}\sum_{|q-p| \leq 2, p\leq Q_c+2}\lambda^{2s+2}_q \int_{\mathbb{R}^3}\Delta_q (u_p\cdot \nabla c_{\leq p-2})c_q\mathrm{d}x\\
&+ \sum_{q \geq -1}\sum_{|q-p| \leq 2, p> Q_c+2}\lambda^{2s+2}_q \int_{\mathbb{R}^3}\Delta_q (u_p\cdot \nabla c_{\leq Q_c})c_q\mathrm{d}x\\
&+ \sum_{q \geq -1}\sum_{|q-p| \leq 2, p> Q_c+2}\lambda^{2s+2}_q \int_{\mathbb{R}^3}\Delta_q (u_p\cdot \nabla c_{(Q_c, p-2]})c_q\mathrm{d}x\\
=: &II_{21}+II_{22}+II_{23}.
\end{split}
\end{equation}
It follows from H\"older's and Young's inequalities that
\begin{equation}\notag
\begin{split}
|II_{21}|\leq & \sum_{q \geq -1}\sum_{|q-p| \leq 2,p\leq Q_c+2}\lambda^{2s+2}_q \int_{\mathbb{R}^3}\big|\Delta_q (u_p\cdot \nabla c_{\leq p-2})c_q\big|\mathrm{d}x\\
\leq & \sum_{q \geq -1}\lambda^{2s+2}_q \|c_q\|_2\sum_{|q-p| \leq 2,p\leq Q_c+2}\| u_p\|_2\| \nabla c_{\leq p-2}\|_\infty\\
\leq &\| \nabla c_{\leq Q_c}\|_\infty \sum_{q \geq -1}\lambda^{s+1}_q \|c_q\|_2 \sum_{|q-p| \leq 2,p\leq Q_c+2}\lambda_p^{s+2}\| u_p\|_2\lambda_p^{-1}\\
\leq & \frac{1}{32} \sum_{q \geq -1}\lambda_q^{2s+4}\|u_q\|_2^2+Cf(t)\sum_{q \geq -1}\lambda_q^{2s+2}\|c_q\|_2^2.
\end{split}
\end{equation}
While the term $II_{22}$ can be treated the same way as $II_{21}$, the third term is estimated as follows, by utilizing the definition of $\Lambda_c(t)$
\begin{equation}\notag
\begin{split}
|II_{23}|\leq &\sum_{q \geq -1}\sum_{|q-p| \leq 2, p> Q_c+2}\lambda^{2s+2}_q \|u_p\|_2\|c_q\|_2 \|\nabla c_{(Q_c, p-2]}\|_\infty\\
\lesssim &\sum_{q >Q_c}\lambda^{2s+2}_q \|u_q\|_2\|c_q\|_2 \|\nabla c_{(Q_c, q]}\|_\infty\\
\lesssim &\sum_{q >Q_c}\lambda^{2s+2}_q \|u_q\|_2\|c_q\|_2 \sum_{Q_c<p\leq q}\lambda_p^{1+\frac3r}\|c_p\|_r\\
\lesssim & C_0\sum_{q >Q_c}\lambda^{2s+2}_q \|u_q\|_2\|c_q\|_2 \sum_{Q_c<p\leq q}\lambda_p\\
\lesssim & C_0\sum_{q >Q_c}\lambda^{s+2}_q \|u_q\|_2\lambda^{s+1}_q\|c_q\|_2 \sum_{Q_c<p\leq q}\lambda_{p-q}\\
\lesssim & C_0\sum_{q >Q_c}\lambda^{2s+4}_q \|u_q\|_2^2+C_0\sum_{q >Q_c}\lambda^{2s+2}_q \|c_q\|_2^2.
\end{split}
\end{equation} 

As for $II_3,$ we first integrate by parts, then split by the wavenumber $Q_u(t)$,
\begin{equation}\notag
\begin{split}
|II_{3}|\leq&\left|\sum_{p\geq -1}\lambda_q^{2s+2}\sum_{-1\leq q\leq p+2}\int_{\mathbb R^3}\Delta_q(u_p\otimes\tilde c_p)\nabla c_q \, \mathrm{d}x\right|\\
\leq&\sum_{p>Q_u}\lambda_q^{2s+2}\sum_{-1\leq q\leq p+2}\int_{\mathbb R^3}\left|\Delta_q(u_p\otimes\tilde c_p)\nabla c_q\right| \, \mathrm{d}x\\
&+\sum_{-1\leq p\leq Q_u}\lambda_q^{2s+2}\sum_{-1\leq q\leq p+2}\int_{\mathbb R^3}\left|\Delta_q(u_p\otimes\tilde c_p)\nabla c_q\right| \, \mathrm{d}x\\
=:& II_{31}+II_{32}.
\end{split}
\end{equation}
Proceeding with the help of H\"older's, Young's and Jensen's inequalities, we have for $s>-\frac12$
\begin{equation}\notag
\begin{split}
|II_{31}|
\lesssim &\sum_{p> Q_u}\|u_p\|_\infty\|c_p\|_2\sum_{-1\leq q\leq p+2}\lambda_q^{2s+3}\|c_q\|_2\\
\lesssim & C_0 \sum_{p> Q_u}\lambda_p\|c_p\|_2\sum_{-1\leq q\leq p+2}\lambda_q^{2s+3}\|c_q\|_2\\
\lesssim & C_0 \sum_{p> Q_u}\lambda_p^{s+2}\|c_p\|_2\sum_{-1\leq q\leq p+2}\lambda_q^{s+2}\|c_q\|_2\lambda_{q-p}^{s+1}\\
\lesssim & C_0 \sum_{p> Q_u}\left(\lambda_p^{2s+4}\|c_p\|_2^2+\left(\sum_{-1\leq q\leq p+2}\lambda_q^{s+2}\|c_q\|_2\lambda_{q-p}^{s+1}\right)^2\right)\\
\lesssim & C_0 \sum_{q\geq -1}\lambda_q^{2s+4}\|c_q\|_2^2,
\end{split}
\end{equation}
and
\begin{equation}\notag
\begin{split}
|II_{32}|
\lesssim &\sum_{-1\leq p\leq Q_u}\|u_p\|_\infty\|c_p\|_2\sum_{-1\leq q\leq p+2}\lambda_q^{2s+3}\|c_q\|_2\\
\lesssim &f(t)\sum_{-1\leq p\leq Q_u}\lambda_p^{-1}\|c_p\|_2\sum_{-1\leq q\leq p+2}\lambda_q^{2s+3}\|c_q\|_2\\
\lesssim &f(t)\sum_{-1\leq p\leq Q_u}\lambda_p^{s+1}\|c_p\|_2\sum_{-1\leq q\leq p+2}\lambda_q^{s+1}\|c_q\|_2\lambda_{q-p}^{s+2}\\
\lesssim &f(t)\sum_{ q \geq -1}\lambda_q^{2s+2}\|c_q\|_2^2.
\end{split}
\end{equation}

We combine the above estimates to conclude that
\begin{equation}\notag
|II|\lesssim C_0 \sum_{q\geq -1}(\lambda_q^{2s+4}\|c_q\|_2^2+\lambda_q^{2s+4}\|u_q\|_2^2)+Q_uf(t)\sum_{q\geq-1}\lambda_q^{2s+2}\|c_q\|_2^2.
\end{equation}
\cbdu

\subsection{Estimate of $III$} 

\begin{Lemma}\label{est-3}
Let $s<0$. We have
\begin{equation}\notag
\begin{split}
|III|\lesssim & C_0 \sum_{q\geq-1}\lambda_q^{2s+2}\|n_q\|_2^2+Q_uf(t)\sum_{q\geq-1}\lambda_q^{2s}\|n_q\|_2^2.
\end{split}
\end{equation}
\end{Lemma}
\pf
We start with Bony's paraproduct decomposition for both of the terms.
\begin{equation}\notag
\begin{split}
III=&\sum_{q\geq -1}\sum_{|q-p|\leq 2}\lambda_q^{2s}\int_{\R^3}\Delta_q(u_{\leq p-2}\cdot\nabla n_p)n_q\, \mathrm{d}x\\
&+\sum_{q\geq -1}\sum_{|q-p|\leq 2}\lambda_q^{2s}\int_{\R^3}\Delta_q(u_{p}\cdot\nabla n_{\leq{p-2}})n_q\, \mathrm{d}x\\
&+\sum_{q\geq -1}\sum_{p\geq q-2}\lambda_q^{2s}\int_{\R^3}\Delta_q(u_p\cdot\nabla\tilde n_p)n_q\, \mathrm{d}x\\
=&III_{1}+III_{2}+III_{3}.
\end{split}
\end{equation}
Since we can estimate $III_1$ and $III_3$ in the same manner as  $II_1$ and $II_3,$ respectively, the details of computation are omitted for simplicity. We claim
\begin{equation}\notag
|III_1|+|III_3|\lesssim  C_0 \sum_{q\geq-1}\lambda_q^{2s+2}\|n_q\|_2^2+Q_uf(t)\sum_{q\geq-1}\lambda_q^{2s}\|n_q\|_2^2.
\end{equation}
We are then left with the estimation of $III_2$ to complete the conclusion. Splitting the term by the wavenumber $Q_u(t),$ we have
\begin{equation}\notag
\begin{split}
III_{2}=& \sum_{ -1 \leq p \leq Q_u}\sum_{|q-p| \leq 2}\lambda^{2s}_q \int_{\mathbb{R}^3}\Delta_q (u_p\cdot \nabla n_{\leq p-2})n_q\mathrm{d}x\\
&+\sum_{ p > Q_u}\sum_{|q-p| \leq 2}\lambda^{2s}_q \int_{\mathbb{R}^3}\Delta_q (u_p\cdot \nabla n_{\leq p-2})n_q\mathrm{d}x\\
=:& III_{21}+III_{22}.
\end{split}
\end{equation}
Applying H\"older's, Young's and Jensen's inequalities, we have for $s<1$
\begin{equation}\notag
\begin{split}
|III_{21}| \leq & \sum_{ -1 \leq p \leq Q_u}\|u_p\|_\infty \sum_{|q-p| \leq 2}\lambda^{2s}_q \|n_q\|_2 \sum_{p' \leq p-2}\lambda_{p'}\|n_{p'}\|_2 \\
\leq & \sum_{ -1 \leq p \leq Q_u}\lambda_p\|u_p\|_\infty \sum_{|q-p| \leq 2}\lambda^{2s-1}_q \|n_q\|_2 \sum_{p' \leq p-2}\lambda_{p'}\|n_{p'}\|_2 \\
\leq & \sum_{ -1 \leq p \leq Q_u}\lambda_p\|u_p\|_\infty \sum_{|q-p| \leq 2}\lambda^{s}_q \|n_q\|_2 \sum_{p' \leq p-2}\lambda_{p'}^{s}\|n_{p'}\|_2\lambda_{q-p'}^{s-1} \\
\lesssim & f(t)\sum_{-1 \leq p \leq Q_u}\Big(\lambda_p^{2s}\|n_q\|_2^2 + \big(\sum_{p' \leq p-2}\lambda_{p'}^{s}\|n_{p'}\|_2\lambda_{p-p'}^{s-1}\big)^2\Big)\\
\lesssim & f(t) \sum_{q \geq -1}\lambda_q^{2s}\|n_q\|_2^2,
\end{split}
\end{equation}
and for $s<0$
\begin{equation}\notag
\begin{split}
|III_{22}| \leq & \sum_{p > Q_u}\|u_p\|_\infty \sum_{|q-p| \leq 2}\lambda^{2s}_q \|n_q\|_2 \sum_{p' \leq p-2}\lambda_{p'}\|n_{p'}\|_2 \\
\leq & \sum_{p > Q_u}\lambda_p^{-1}\|u_p\|_\infty \sum_{|q-p| \leq 2}\lambda^{2s+1}_q \|n_q\|_2 \sum_{p' \leq p-2}\lambda_{p'}\|n_{p'}\|_2 \\
\leq & C_0 \sum_{q>Q_u-2}\lambda^{s+1}_q \|n_q\|_2 \sum_{p' \leq q}\lambda_{p'}^{s+1}\|n_{p'}\|_2\lambda_{q-p'}^{s} \\
\lesssim & C_0 \sum_{q \geq -1}\lambda_q^{2s+2}\|n_q\|_2^2.
\end{split}
\end{equation}

\cbdu

\subsection{Estimate of $IV$}
Here we use H\"older's and Young's inequalities to obtain
\begin{equation}\label{estnn}
\begin{split}
|IV| \leq & \sum_{q \geq -1}\lambda_q^{s+1}\|n_q\|_2\lambda_q^{s+1}\|u_q\|_2\\
\leq & \frac1{32} \sum_{q \geq -1}\lambda_q^{2s+2}\|n_q\|_2^2+ C \sum_{q \geq -1}\lambda_q^{2s+2}\|u_q\|_2^2
\end{split}
\end{equation}
for an absolute constant $C$.

\subsection{Estimate of $V$}
\begin{Lemma}\label{ppsnc}
If $r> 3$ and $s<0$, we have 
\begin{equation}\notag
\begin{split}
|V| \leq &\left(\frac1{32}+CC_0\right)\sum_{q\geq -1}\left(\lambda_q^{2s+4}\|c_q\|_2^2+\lambda_q^{2s+2}\|n_q\|_2^2\right)\\
&+C(f(t)+1)\sum_{q\geq-1}\left(\lambda_q^{2s+2}\|c_q\|_2^2+\lambda_q^{2s}\|n_q\|_2^2\right).
\end{split}
\end{equation}
\end{Lemma} 
\pf
Bony's paraproduct decomposition leads to
\begin{equation}\notag
\begin{split}
V=& \sum_{q \geq -1} \sum_{|q-p| \leq 2}\lambda^{2s+2}_q\int_{\mathbb{R}^3}\Delta_q(n_{\leq p-2}c_p)c_q\mathrm{d}x\\
&+ \sum_{q \geq -1} \sum_{|q-p| \leq 2}\lambda^{2s+2}_q\int_{\mathbb{R}^3}\Delta_q(n_pc_{\leq p-2})c_q\mathrm{d}x\\
&+ \sum_{q \geq -1} \sum_{p \geq q-2}\lambda^{2s+2}_q\int_{\mathbb{R}^3}\Delta_q(\tilde n_pc_p)c_q\mathrm{d}x\\
=: & V_1+V_2+V_3.
\end{split}
\end{equation}
We further split $V_1$ by the wavenumber $\Lambda_c(t)$
\begin{equation}\notag
\begin{split}
V_1=& \sum_{q \leq Q_c} \sum_{|q-p| \leq 2}\lambda^{2s+2}_q\int_{\mathbb{R}^3}\Delta_q(n_{\leq p-2}c_p)c_q\mathrm{d}x\\
& +\sum_{q >Q_c} \sum_{|q-p| \leq 2}\lambda^{2s+2}_q\int_{\mathbb{R}^3}\Delta_q(n_{\leq p-2}c_p)c_q\mathrm{d}x\\
=:&V_{11}+V_{12}.
\end{split}
\end{equation}
To estimate $V_{11}$, H\"older's inequality gives 
\begin{equation}\notag
\begin{split}
|V_{11}| \leq & \sum_{q \leq Q_c} \lambda^{2s+2}_q\|c_q\|_\infty\sum_{|q-p| \leq 2}\|c_p\|_2\sum_{p' \leq p-2}\|n_{p'}\|_2 \\
\lesssim & \sum_{q \leq Q_c}\lambda_q\|c_q\|_\infty\sum_{|q-p| \leq 2}\lambda_p^{s+2}\|c_p\|_2\sum_{p' \leq p-2}\lambda_{p'}^{s} \|n_{p'}\|_2\lambda_{p-p'}^{s-1}\lambda_{p'}^{-1}.\\
\end{split}
\end{equation}
Since $s-1<0$ for $s<0$, we apply Young's and Jensen's inequalities to obtain 
\begin{equation}\notag
\begin{split}
|V_{11}|
\leq & \frac1{32} \sum_{q \leq Q_c+2}\lambda_q^{2s+4}\|c_q\|_2^2+C\sum_{q \leq Q_c+2}\lambda_q^{2}\|c_q\|_\infty^2\left(\sum_{p' \leq q}\lambda_{p'}^{s}\|n_{p'}\|_2\lambda_{q-p'}^{s-1}\right)^2\\
\leq & \frac1{32} \sum_{q \leq Q_c+2}\lambda_q^{2s+4}\|c_q\|_2^2+Cf(t)\sum_{q \leq Q_c+2}\sum_{p' \leq q}\lambda_{p'}^{2s}\|n_{p'}\|_2^2\lambda_{q-p'}^{s-1}\\
\leq & \frac1{32} \sum_{q \leq Q_c+2} \lambda_q^{2s+4}\|c_q\|_2^2+Cf(t)\sum_{q \leq Q_c}\lambda_q^{2s}\|n_q\|_2^2.
\end{split}
\end{equation}
Regarding $V_{12}$, we have 
\begin{equation}\notag
\begin{split}
|V_{12}| \leq & \sum_{q > Q_c} \lambda^{2s+2}_q\|c_q\|_\infty\sum_{|q-p| \leq 2}\|c_p\|_2\sum_{p' \leq p-2}\|n_{p'}\|_2 \\
 \lesssim & \sum_{q > Q_c} \lambda^{2s+2+\frac3r}_q\|c_q\|_r\sum_{|q-p| \leq 2}\|c_p\|_2\sum_{p' \leq p-2}\|n_{p'}\|_2 \\
\lesssim &C_0 \sum_{q > Q_c-2} \lambda^{2s+2}_q\|c_q\|_2\sum_{p' \leq q}\|n_{p'}\|_2 \\
\lesssim &C_0 \sum_{q > Q_c-2} \lambda^{s+2}_q\|c_q\|_2\sum_{p' \leq q}\lambda_{p'}^s\|n_{p'}\|_2 \lambda_{q-p'}^{s}.
\end{split}
\end{equation}
Again, since $s<0$, we can apply Young's and Jensen's inequalities
\begin{equation}\notag
\begin{split}
|V_{12}| 
\lesssim &C_0 \sum_{q > Q_c-2} \lambda^{2s+4}_q\|c_q\|_2^2+\sum_{q>Q_c}\left(\sum_{p' \leq q}\lambda_{p'}^s\|n_{p'}\|_2 \lambda_{q-p'}^{s}\right)^2\\
\lesssim &C_0 \sum_{q > Q_c-2} \lambda^{2s+4}_q\|c_q\|_2^2+\sum_{q\geq-1}\lambda_{p'}^{2s}\|n_{p'}\|_2^2.
\end{split}
\end{equation}
We also split $V_2$ by the wavenumber $\Lambda_c(t)$ as
\begin{equation}\notag
\begin{split}
V_2
=&\sum_{q \geq -1} \sum_{|q-p| \leq 2}\lambda^{2s+2}_q\int_{\mathbb{R}^3}\Delta_q(n_pc_{\leq p-2})c_q\mathrm{d}x\\
=&\sum_{q \geq -1} \sum_{|q-p| \leq 2,p\leq Q_c+2}\lambda^{2s+2}_q\int_{\mathbb{R}^3}\Delta_q(n_pc_{\leq p-2})c_q\mathrm{d}x\\
&+\sum_{q \geq -1} \sum_{|q-p| \leq 2,p> Q_c+2}\lambda^{2s+2}_q\int_{\mathbb{R}^3}\Delta_q(n_pc_{\leq Q_c})c_q\mathrm{d}x\\
&+\sum_{q \geq -1} \sum_{|q-p| \leq 2,p> Q_c+2}\lambda^{2s+2}_q\int_{\mathbb{R}^3}\Delta_q(n_pc_{(Q_c,p-2]})c_q\mathrm{d}x\\
=:&V_{21}+V_{22}+V_{23}.
\end{split}
\end{equation}
H\"older's and Young's inequalities yield the following estimate on $V_{21}$, 
\begin{equation}\notag
\begin{split}
|V_{21}| \lesssim & \sum_{q \geq -1} \lambda^{2s+2}_q\|c_q\|_2 \sum_{|q-p| \leq 2}\|n_p\|_2\sum_{p' \leq p-2\leq Q_c}\|c_{p'}\|_\infty \\
\lesssim & \sum_{q \geq -1}\lambda_q^{s+1}\|c_q\|_2 \sum_{|q-p| \leq 2}\lambda_p^{s+1} \|n_p\|_2\sum_{p' \leq p-2\leq Q_c}\lambda_{p'}\|c_{p'}\|_\infty\lambda_{p'}^{-1}\\
\leq & \frac1{32}\sum_{q \geq -1}\lambda_q^{2s+2}\|n_q\|_2^2+C\sum_{p\geq-1}\lambda_p^{2s+2}\|c_p\|_2^2\left(\sum_{p'\leq p-2\leq Q_c}\lambda_{p'}\|c_{p'}\|_\infty\lambda_{p'}^{-1}\right)^2\\
\leq & \frac1{32}\sum_{q \geq -1}\lambda_q^{2s+2}\|n_q\|_2^2+C\sum_{p\geq-1}\lambda_p^{2s+2}\|c_p\|_2^2\sum_{p'\leq p-2\leq Q_c}\lambda_{p'}^2\|c_{p'}\|_\infty^2\lambda_{p'}^{-1}\\
\leq & \frac1{32}\sum_{q \geq -1}\lambda_q^{2s+2}\|n_q\|_2^2+Cf(t)\sum_{p\geq-1}\lambda_p^{2s+2}\|c_p\|_2^2\sum_{p'\leq p-2\leq Q_c}\lambda_{p'}^{-1}\\
\leq & \frac1{32}\sum_{q \geq -1} \lambda_q^{2s+2}\|n_q\|_2^2+Cf(t)\sum_{q \geq -1}\lambda_q^{2s+2}\|c_q\|_2^2.
\end{split}
\end{equation}
Note that $V_{22}$ can be estimated in the same way. Provided $\frac3r\leq 1$, the last term $V_{23}$ is estimated as follows,
\begin{equation}\notag
\begin{split}
|V_{23}|\leq & \sum_{q \geq -1} \sum_{|q-p| \leq 2,p> Q_c+2}\lambda^{2s+2}_q \|n_p\|_2\|c_q\|_{\frac{2r}{r-2}}\|c_{(Q_c,p-2]}\|_r\\
\lesssim & \sum_{q \geq -1} \lambda^{2s+2+\frac3r}_q \|n_q\|_2\|c_q\|_2\sum_{q> Q_c,Q_c<p'\leq q}\|c_{p'}\|_r\\
\lesssim & C_0\sum_{q \geq -1} \lambda^{2s+2+\frac3r}_q \|n_q\|_2\|c_q\|_2\sum_{q> Q_c,Q_c<p'\leq q}\lambda_{p'}^{-\frac3r}\\
\lesssim & C_0\sum_{q \geq -1} \lambda^{s+1}_q \|n_q\|_2\lambda^{s+2}_q\|c_q\|_2\lambda_q^{\frac3r-1}\sum_{q> Q_c,Q_c<p'\leq q}\lambda_{p'}^{-\frac3r}\\
\lesssim & C_0\sum_{q \geq -1} \lambda^{s+1}_q \|n_q\|_2\lambda^{s+2}_q\|c_q\|_2\\
\lesssim & C_0\sum_{q \geq -1} \lambda^{2s+2}_q \|n_q\|_2^2+C_0\sum_{q \geq -1}\lambda^{2s+4}_q\|c_q\|_2^2.\\
\end{split}
\end{equation}

The term $V_3$ can be dealt with in a similar way as $V_1$. We first split the sum,
\begin{equation}\notag
\begin{split}
V_3=&\sum_{q \leq Q_c} \sum_{p \geq q-2}\lambda^{2s+2}_q\int_{\mathbb{R}^3}\Delta_q(\tilde n_pc_p)c_q\mathrm{d}x\\
&+\sum_{q < Q_c} \sum_{p \geq q-2}\lambda^{2s+2}_q\int_{\mathbb{R}^3}\Delta_q(\tilde n_pc_p)c_q\mathrm{d}x.\\
\end{split}
\end{equation}
Without giving details, we claim that 
\begin{equation}\notag
|V_3|\leq \left(\frac1{32}+CC_0\right)\sum_{q\geq -1}\lambda_q^{2s+4}\|c_q\|_2^2+C(f(t)+1)\sum_{q\geq-1}\lambda_q^{2s}\|n_q\|_2^2.
\end{equation}

\cbdu

\subsection{Estimate of $VI$}
\begin{Lemma}\label{ppscc} Let $-\frac12<s<0$ and $3<r<\frac 3{1+s}$. We have
\begin{equation}\notag
|VI| \leq \left(\frac1{32}+CC_0\right)\sum_{q\geq -1}\lambda_q^{2s+2}\|n_q\|_2^2+Cf(t)\sum_{q\geq -1}\lambda_q^{2s}\|n_q\|_2^2.
\end{equation}
\end{Lemma}
\pf
Utilizing Bony's paraproduct, $VI$ can be decomposed as
\begin{equation}\notag
\begin{split}
VI=
&-\sum_{q\geq -1}\sum_{|q-p|\leq 2}\lambda_q^{2s}\int_{\R^3}\Delta_q(n_{\leq p-2}\nabla c_p)\nabla n_q\, \mathrm{d}x\\
&-\sum_{q\geq -1}\sum_{|q-p|\leq 2}\lambda_q^{2s}\int_{\R^3}\Delta_q(n_{p}\nabla c_{\leq{p-2}})\nabla n_q\, \mathrm{d}x\\
&-\sum_{q\geq -1}\sum_{p\geq q-2}\lambda_q^{2s}\int_{\R^3}\Delta_q(\tilde n_p\nabla c_p)\nabla n_q\, \mathrm{d}x\\
=:& VI_{1}+VI_{2}+VI_{3}.
\end{split}
\end{equation}
We continue to decompose $VI_1$ by $Q_c$,
\begin{equation}\notag
\begin{split}
VI_1=&-\sum_{q\geq -1}\sum_{|q-p|\leq 2, p\leq Q_c}\lambda_q^{2s}\int_{\R^3}\Delta_q(n_{\leq p-2}\nabla c_p)\nabla n_q\, \mathrm{d}x\\
&-\sum_{q\geq -1}\sum_{|q-p|\leq 2, p> Q_c}\lambda_q^{2s}\int_{\R^3}\Delta_q(n_{\leq p-2}\nabla c_p)\nabla n_q\, \mathrm{d}x\\
=:&VI_{11}+VI_{12}.
\end{split}
\end{equation}
To estimate $VI_{11}$, we apply H\"older's inequality first, 
\begin{equation}\notag
\begin{split}
|VI_{11}| \lesssim& \sum_{q \geq -1}\lambda_q^{2s+1}\|n_q\|_2 \sum_{|q-p| \leq 2,p\leq Q_c}\|\nabla c_p\|_\infty \sum_{p' \leq p-2}\|n_{p'}\|_2\\
\lesssim& \sum_{q \geq -1}\lambda_q^{s+1}\|n_q\|_2 \sum_{|q-p| \leq 2,p\leq Q_c}\|\nabla c_p\|_\infty\sum_{p' \leq p-2}\lambda_{p'}^{s}\|n_{p'}\|_2\lambda_{p-p'}^{s}.\\
\end{split}
\end{equation}
We then apply Young's and Jensen's inequalities,
\begin{equation}\notag
\begin{split}
|VI_1| \leq & \frac1{32} \sum_{q \geq -1}\lambda_q^{2s+2}\|n_q\|_2^2+ Cf(t)\sum_{q \geq -1}\big(\sum_{p' \leq q}\lambda_{p'}^{s}\|n_{p'}\|_2\lambda_{q-p'}^{s}\big)^2\\
\leq & \frac1{32} \sum_{q \geq -1}\lambda_q^{2s+2}\|n_q\|_2^2+Cf(t)\sum_{q \geq -1}\lambda_q^{2s}\|n_q\|_2^2,
\end{split}
\end{equation}
where  we require $s <0$ to ensure $\sum_{p'\leq p-2} \lambda^{s}_{p-p'} <\infty$.

Again, applying H\"older's inequality first to $VI_{12}$ yields
\begin{equation}\notag
\begin{split}
|VI_{12}|\leq &\sum_{q\geq -1}\sum_{|q-p|\leq 2, p> Q_c}\lambda_q^{2s}\|\nabla n_q\|_2\|\nabla c_p\|_r\|n_{\leq p-2}\|_{\frac{2r}{r-2}}\\
\lesssim& \sum_{q\geq -1}\sum_{|q-p|\leq 2, p> Q_c}\lambda_q^{2s+2-\frac3r}\|n_q\|_2\lambda_p^{\frac3r}\|c_p\|_r\sum_{p'\leq p-2}\lambda_{p'}^{\frac3r}\|n_{p'}\|_2\\
\lesssim& C_0\sum_{q\geq -1}\lambda_q^{2s+2-\frac3r}\|n_q\|_2\sum_{p'\leq q}\lambda_{p'}^{\frac3r}\|n_{p'}\|_2\\
\lesssim& C_0\sum_{q\geq -1}\lambda_q^{s+1}\|n_q\|_2\sum_{p'\leq q}\lambda_{p'}^{s+1}\|n_{p'}\|_2\lambda_{q-p}^{s+1-\frac3r}.\\
\end{split}
\end{equation}
Following Young's and Jensen's inequalities, we obtain for $s+1<\frac3r$
\begin{equation}\notag
\begin{split}
|VI_{12}|\leq &C C_0\sum_{q\geq -1}\lambda_q^{2s+2}\|n_q\|_2^2+CC_0\sum_{q\geq-1}\left(\sum_{p'\leq q}\lambda_{p'}^{s+1}\|n_{p'}\|_2\lambda_{q-p}^{s+1-\frac3r}\right)^2\\
\leq &C C_0\sum_{q\geq -1}\lambda_q^{2s+2}\|n_q\|_2^2+CC_0\sum_{q\geq-1}\sum_{p'\leq q}\lambda_{p'}^{2s+2}\|n_{p'}\|_2^2\lambda_{q-p}^{s+1-\frac3r}\\
\leq &C C_0\sum_{q\geq -1}\lambda_q^{2s+2}\|n_q\|_2^2.
\end{split}
\end{equation}

The first step of dealing with $VI_2$ is to split it by $Q_c$ as well,
\begin{equation}\notag
\begin{split}
VI_2=&-\sum_{q\geq -1}\sum_{|q-p|\leq 2,p\leq Q_c+2}\lambda_q^{2s}\int_{\R^3}\Delta_q(n_{p}\nabla c_{\leq{p-2}})\nabla n_q\, \mathrm{d}x\\
&-\sum_{q\geq -1}\sum_{|q-p|\leq 2,p> Q_c+2}\lambda_q^{2s}\int_{\R^3}\Delta_q(n_{p}\nabla c_{\leq Q_c})\nabla n_q\, \mathrm{d}x\\
&-\sum_{q\geq -1}\sum_{|q-p|\leq 2,p> Q_c+2}\lambda_q^{2s}\int_{\R^3}\Delta_q(n_{p}\nabla c_{(Q_c,p-2]})\nabla n_q\, \mathrm{d}x\\
=:&VI_{21}+VI_{22}+VI_{23}.
\end{split}
\end{equation}
Using H\"older's and Young's inequalities, and the fact that $\lambda_{q-p}^{s} \sim C$ as $p,q$ are close to each other, we have 
\begin{equation}\notag
\begin{split}
|VI_{21}| \leq& \sum_{q \geq -1}\lambda_q^{2s+1}\|n_q\|_2 \sum_{|q-p| \leq 2,p\leq Q_c+2}\| n_p\|_2 \|\nabla c_{\leq p-2}\|_\infty\\
\leq& \sum_{q \geq -1}\lambda_q^{s+1}\|n_q\|_2 \sum_{|q-p| \leq 2,p\leq Q_c+2}\lambda_p^{s}\| n_p\|_2 \|\nabla c_{\leq p-2}\|_\infty\lambda_{q-p}^{s}\\
\leq & \frac1{32} \sum_{q \geq -1}\lambda_q^{2s+2}\|n_q\|_2^2+Cf(t)\sum_{q \geq -1}\lambda_q^{2s}\|n_q\|_2^2.
\end{split}
\end{equation}
We skip the computation for $VI_{22}$ which can be estimated in a likely way as $VI_{21}$. We proceed to estimate $V_{23}$, provided $r>3$
\begin{equation}\notag
\begin{split}
|VI_{23}|\leq &\sum_{q\geq -1}\sum_{|q-p|\leq 2,p> Q_c+2}\lambda_q^{2s}\|n_p\|_2\|\nabla c_{(Q_c,p-2]}\|_{r}\|\nabla n_q\|_{\frac{2r}{r-2}}\\
\lesssim &\sum_{q\geq -1}\lambda_q^{2s+1+\frac3r}\|n_q\|_2^2\sum_{ Q_c<p'\leq q}\lambda_{p'}\|c_{p'}\|_{r}\\
\lesssim &C_0\sum_{q\geq -1}\lambda_q^{2s+1+\frac3r}\|n_q\|_2^2\sum_{ Q_c<p'\leq q}\lambda_{p'}^{1-\frac3r}\\
\lesssim &C_0\sum_{q\geq -1}\lambda_q^{2s+2}\|n_q\|_2^2\sum_{ Q_c<p'\leq q}\lambda_{p'-q}^{1-\frac3r}\\
\lesssim &C_0\sum_{q\geq -1}\lambda_q^{2s+2}\|n_q\|_2^2.
\end{split}
\end{equation}

Finally, we notice that $VI_3$ can be handled in an analogous way of $VI_1$; thus the details of computation are omitted. It completes the proof of the lemma.


\cbdu 

Summing inequalities in Lemma \ref{ppsuu}- Lemma \ref{ppscc}  and (\ref{eqtnn})-(\ref{estnn}) produces the following Gr\"onwall type of inequality,
\begin{equation}\notag
\begin{split}
&\frac{\mathrm{d}}{\mathrm{d}t}\sum_{q \geq -1}(\lambda_q^{2s}\|n_q\|_2^2+\lambda_q^{2s+2}\|c_q\|_2^2+\lambda_q^{2s+2}\|u_q\|_2^2)\\
\leq & (-2+CC_0) \sum_{q\geq -1}(\lambda_q^{2s+2}\|n_q\|_2^2+\lambda_q^{2s+4}\|c_q\|_2^2+\lambda_q^{2s+4}\|u_q\|_2^2)\\
&+CQ_uf(t)\sum_{q\geq -1}(\lambda_q^{2s}\|n_q\|_2^2+\lambda_q^{2s+2}\|c_q\|_2^2+\lambda_q^{2s+2}\|u_q\|_2^2).
\end{split}
\end{equation}
Thus, $C_0$ can be chosen small enough such that $CC_0<\frac14$. On the other hand, combining the definition of $\Lambda_u(t)$ and Bernstein's inequality, one can deduce
\[1\leq C_0^{-1}\Lambda_u^{-1}\|u_{Q_u}\|_\infty\lesssim C_0^{-1}\Lambda_u^{\frac12}\|u_{Q_u}\|_2\lesssim C_0^{-1}\Lambda_u^{-\frac12-s}\|u\|_{\dot H^{s+1}},\]
which implies that for $s>-\frac12$,
\[Q_u=\log \Lambda_u\lesssim 1+\log {\|u\|_{\dot H^{s+1}}}.\]
Hence the energy inequality becomes 
\begin{equation}\notag
\begin{split}
&\frac{\mathrm{d}}{\mathrm{d}t}\sum_{q \geq -1}(\lambda_q^{2s}\|n_q\|_2^2+\lambda_q^{2s+2}\|c_q\|_2^2+\lambda_q^{2s+2}\|u_q\|_2^2)\\
\leq &-\sum_{q\geq -1}(\lambda_q^{2s+2}\|n_q\|_2^2+\lambda_q^{2s+4}\|c_q\|_2^2+\lambda_q^{2s+4}\|u_q\|_2^2)\\
 &+C f(t)\left(1+\log {\|u\|_{\dot H^{s+1}}}\right)\sum_{q\geq -1}(\lambda_q^{2s}\|n_q\|_2^2+\lambda_q^{2s+2}\|c_q\|_2^2+\lambda_q^{2s+2}\|u_q\|_2^2).
\end{split}
\end{equation}
By the hypothesis (\ref{c1u1q}) and Gr\"onwall's inequality, we can conclude that 
\begin{align*}
n\in L^\infty(0,T; H^{s})\cap L^2(0,T; H^{s+1}),\\
c\in L^\infty(0,T; H^{s+1})\cap L^2(0,T; H^{s+2}),\\
u\in L^\infty(0,T; H^{s+1})\cap L^2(0,T; H^{s+2}).
\end{align*}
We consider a particular case $s=-\varepsilon$ for small enough $\varepsilon>0$. We realize that our solution is regular via bootstrapping arguments. In fact, we have
\[n\in L^\infty(0,T; H^{-\varepsilon}), \ \ u\in L^\infty(0,T; H^{1-\varepsilon}), \ \ c\in L^\infty(0,T; H^{1-\varepsilon})\cap L^\infty(0,T; L^\infty).\]
By scaling, it is known that $H^{-\frac 12}$, $H^{\frac 12}$, and $H^{\frac 32}$ are critical for $n, u$ and $c$, respectively. For small enough $\varepsilon>0$, $H^{-\varepsilon}$ is subcritical for $n$; so is $H^{1-\varepsilon}$ for $u$. Thus, it suffices to bootstrap the equation of $c$ to obtain higher regularity for $c$. We recall
\[c_t-\Delta c=-u\cdot\nabla c-nc.\]
Since $u\in L^\infty(0,T; H^{1-\varepsilon})$ and $\nabla c\in L^2(0,T; H^{1-\varepsilon})$,  by Sobolev embedding $H^{1-\varepsilon} \hookrightarrow L^{\frac{6}{1+2\varepsilon}},$
we have
$u\cdot\nabla c\in L^2(0,T; L^{\frac{3}{1+2\varepsilon}})$. Similarly, the fact of $c\in L^\infty(0,T; H^{1-\varepsilon})$ and $n\in L^2(0,T; H^{1-\varepsilon})$ implies $nc\in L^2(0,T; L^{\frac{3}{1+2\varepsilon}})$. Then the standard maximal regularity theory of heat equation yields 
\[c\in H^1(0,T; L^{\frac{3}{1+2\varepsilon}})\cap L^2(0,T; W^{2, \frac{3}{1+2\varepsilon}}).\]
We need to bootstrap one more time. Now we have $\nabla c\in L^2(0,T; W^{1, \frac{3}{1+2\varepsilon}})$ which, along with $u\in L^\infty(0,T; H^{1-\varepsilon})$, implies $u\cdot\nabla c\in L^2(0,T; L^{\frac{6}{1+6\varepsilon}})$. On the other hand, we have $nc\in L^2(0,T; L^{\frac{6}{1+2\varepsilon}})$ from the estimate $n\in L^2(0,T; H^{1-\varepsilon})$ and the maximal principle $c\in L^\infty(0,T;L^\infty)$.  Again, the maximal regularity theory of heat equation produces that 
\[c\in H^1(0,T; L^{\frac{6}{1+6\varepsilon}})\cap L^2(0,T; W^{2, \frac{6}{1+6\varepsilon}}).\]
As a consequence of the mixed derivative theorem \cite{PS}, we have
\[c\in W^{1-\theta, 2}(0,T; W^{2\theta,  \frac{6}{1+6\varepsilon}})\]
for any $\theta\in[0,1]$. In fact, if we take $\theta\in (\frac{1+6\varepsilon}4,\frac12)$, Sobolev embedding theorem shows that
\[c\in W^{1-\theta, 2}(0,T; W^{2\theta,  \frac{6}{1+6\varepsilon}})\hookrightarrow L^\infty(0,T; H^{\frac32+\varepsilon_0})\]
for an small enough constant $\varepsilon_0>0$. Notice that $H^{\frac 32}$ is critical for $c$. Thus we can stop the bootstrapping for $c$ equation. Regarding the density function $n$, although the obtained estimates are in subcritical space already, we would like to further improve the estimates to reach spaces with even higher regularity (i.e. Sobolev spaces with positive smoothness index).
Indeed, we look at the $n$ equation again,
\[n_t-\Delta n=-u\cdot \nabla n-\nabla\cdot (n\nabla c).\]
Due to the fact $\nabla c\in L^\infty(0,T; H^{\frac12+\varepsilon_0})$ and $n\in L^2(0,T; H^{1-\varepsilon})$, Sobolev embedding theorem yields $n\nabla c\in L^2(0,T; L^2)$ and hence $\nabla \cdot (n\nabla c) \in L^2(0,T; H^{-1})$. While $u\in L^\infty(0,T; H^{1-\varepsilon})$ and $n\in L^2(0,T; H^{1-\varepsilon})$ together imply $un\in L^2(0,T; L^{\frac3{1+2\varepsilon}})$ which is embedded in $L^2(0,T; L^2)$.  Thus $u\cdot \nabla n=\nabla\cdot (un)\in L^2(0,T; H^{-1})$. Applying Lemma \ref{le-heat} with $\alpha=-1$ to the $n$ equation, we claim 
\[n\in H^1(0,T; H^{-1})\cap L^2(0,T; H^1).\] 
Summarizing the analysis above gives us
\begin{equation}\notag
\begin{split}
&n\in L^\infty(0,T; H^{-\varepsilon})\cap L^2(0,T; H^{1})\cap H^1(0,T; H^{-1}),\\
&u\in L^\infty(0,T; H^{1-\varepsilon})\cap L^2(0,T; H^{2-\varepsilon}),\\
&c\in L^\infty(0,T; H^{\frac32+\varepsilon_0})\cap L^2(0,T; W^{2, \frac{6}{1+6\varepsilon}}).
\end{split}
\end{equation}
Since each of the three functions $n,u$ and $c$ is in higher regularity space than its critical Sobolev space, further bootstrapping procedures for parabolic equations and standard argument of extending regularity can be applied to infer that the solution 
$(u,n,c)$ is regular up to time $T$.

\bigskip

{\textbf{Acknowledgement. }}
The authors would like to thank Professor Gieri Simonett for pointing out a mistake in the early version of the manuscript.



\end{document}